\documentclass{amsart}
\usepackage{amssymb}
\newtheorem{thm}{Theorem}
\newtheorem{lem}{Lemma}
\newtheorem{prop}{Proposition}
\newtheorem{cor}{Corollary}

\begin{document}

\pagestyle{myheadings} \markboth{\rm Milo\v{s} Arsenovi\'c and Romi F. Shamoyan}{\rm Embedding relations and
boundedness...}

\bibliographystyle{plain}

\title[Embedding relations and boundedness...]{Embedding relations and boundedness of the multifunctional operators in tube domains over symmetric cones}

\author{Romi F. Shamoyan and Milo\v{s} Arsenovi\'{c}}

\thanks{Supported by Ministry of Science, Serbia, project M144010}

\date{}

\maketitle

\begin{abstract}

We obtain a new sufficient condition for the continuity of the Bergman projection in tube domains over
symmetric cones using multifunctional embeddings. We also obtain some embedding relations between the
generalized Hilbert-Hardy spaces and the mixed-norm Bergman spaces in this setting.

\end{abstract}

\footnotetext[1]{Mathematics Subject Classification 2010 Primary 42B35 Secondary 32A07, 32M15  Key words
and Phrases: Hardy spaces, Bergman spaces, Symmetric cones, Bergman projection, tube domains}

\section{Introduction and statements of the results}

Let $T_\Omega = V + i\Omega$ be the tube domain over an irreducible symmetric cone $\Omega$ in the
complexification $V^{\mathbb C}$ of an $n$-dimensional euclidean space $V$. Following the notation of
\cite{FaKo} we denote the rank of the cone $\Omega$ by $r$ and by $\Delta$ the determinant function on $V$.
Letting $V = \mathbb R^n$, we have as an example of a symmetric cone on $\mathbb R^n$ the Lorentz cone
$\Lambda_n$ which is a rank $2$ cone defined for $n \geq 3$ by
$$\Lambda_n = \{ y \in \mathbb R^n : y_1^2 - \cdots - y_n^2 > 0, y_1 > 0 \}.$$
The determinant function in this case is given by the Lorentz form
$$\Delta(y) = y_1^2 - \cdots - y_n^2.$$

Let us introduce some convenient notation regarding multi-indices.

If $t = (t_1, \ldots, t_r)$, then $t^\star = (t_r, \ldots, t_1)$ and, for $a \in \mathbb R$, $t+a =
(t_1 + a, \ldots, t_n + a)$. Also, if $t, k \in \mathbb R^n$, then $t < k$ means $t_j < k_j$ for all
$1 \leq j \leq r$.

We are going to use the following multi-index
$$g_0 = \left( (j-1)\frac{d}{2} \right)_{1 \leq j \leq r}, \;\;\; \mbox{\rm where} \;\;\;\
(r-1) \frac{d}{2} = \frac{n}{r} - 1.$$

For $1 \leq p, q < + \infty$ and $\nu \in \mathbb R^r$, we denote by $A_\nu^{p,q}(T_\Omega)$ the mixed-norm
Bergman space consisting of analytic functions $f$ in $T_\Omega$ such that
$$ \| f \|_{L_\nu^{p,q}} = \left( \int_\Omega \left( \int_V |F(x+iy)|^p d\,x \right)^{q/p} \Delta_\nu (y)
\frac{d\, y}{\Delta (y)^{n/r}} \right)^{1/q} < \infty,$$
where $\Delta_\nu$ is the generalized power function to be defined in the next section. The space $A_\nu^{p,q}(T_\Omega)$ is nontrivial if and only if $\nu > g_0$, see \cite{DD}. When $p = q$ we write
$A_\nu^{p,q}(T_\Omega) = A^p_\nu(T_\Omega)$; the classical Bergman space $A^p(\Omega)$ corresponds to
$\nu = (n/r, \ldots, n/r)$.

The (weighted) Bergman projection $P_\nu$ is the orthogonal projection from the Hilbert space $L^2_\nu(T_\Omega)$
onto its closed subspace $A^2_\nu(T_\Omega)$ and it is given by the following integral formula
\begin{equation}\label{bpro}
P_\nu f(z) = d_\nu \int_{T_\Omega} B_\nu (z, w) f(w) d V_\nu (w),
\end{equation}
where $B_\nu (z, w) = c_\nu \Delta^{-(\nu + \frac{n}{r})}((z-\overline w)/i)$ is the Bergman reproducing kernel
for $A^2_\nu$, see \cite{FaKo}. Here we used notation $d V_\nu (w) = \Delta^{\nu - \frac{n}{r}}(v) du dv$,
where $w = u +iv \in T_\Omega$.

The problem of boundedness of the Bergman projection on tube domains over symmetric cones has been considered by
several authors (see \cite{BeBo}, \cite{BBPR}, \cite{BBGR}, \cite{BBGRS} and references therein) and still
remains open. The best known results have been obtained in \cite{GS} in the setting of the light cone. Recently,
an equivalent condition for the boundedness of the Bergman projection in terms of Hardy-type inequalities and duality
was obtained in \cite{BBGRS}. We introduce here the operators $T_\beta$, $\beta = (\beta_1, \ldots, \beta_m)$ which
generalize the Bergman projection and are defined by
$$T_\beta(\overrightarrow f)(\overrightarrow z) = \int_{T_\Omega}
\frac{ \left(\prod_{j=1}^m f_j(z) \right) \Delta^{ \frac{1}{m} \sum_{j=1}^m \beta_j}(\Im z)}
{\prod_{j=1}^m \Delta^{\frac{1}{m}(\frac{n}{r} + \beta_j)}(\frac{z_j - \overline z}{i})}
\frac{d V(z)}{\Delta^{\frac{n}{r}}(\Im z)},$$
where $\overrightarrow f = (f_1, \ldots, f_n)$, $\overrightarrow z = (z_1, \ldots, z_n)$, $z_j \in T_\Omega$ and
$f_j \in L^1_{loc}(T_\Omega)$ for $1 \leq j \leq m$. Combining classical arguments with integrability properties of
the Bergman kernel and determinant function we obtain the following sufficient condition for the boundedness of the
operator $T_\beta$ from the product space
$$\prod_{k=1}^m L^p_{m\nu_k + (m-1)\frac{n}{r}}(T_\Omega) = L^p_{m\nu_1 + (m-1)\frac{n}{r}}(T_\Omega) \times
\cdots \times L^p_{m\nu_m + (m-1)\frac{n}{r}}(T_\Omega)$$
to the space $L^p((T_\Omega)^m, \prod_{k=1}^m \Delta^{\nu_k - \frac{n}{r}} d V(z_k))$. The idea to consider such
multifunctional operator is motivated by \cite{LS}. Some results of this paper are analogous to results of \cite{LS}
proven in the case of the unit ball in $\mathbb C^n$. We note here that almost all multifunctional results of this paper are well known in the case $m=1$. For example, the case $m=1$ of the following theorem is contained in \cite{BBPR}.

\begin{thm}
Let $\nu_k \in \mathbb R$, $k=1,\ldots, m$, $m>1$, $1 \leq p < \infty$ and $\beta = (\beta_1,\ldots, \beta_n)$. If the
parameters satisfy the following conditions
\begin{equation}\label{cond1}
\frac{1}{m} \sum_{j=1}^m \beta_j > \frac{n}{r} - 1,
\end{equation}

\begin{equation}\label{cond2}
1 \leq p < 1 + m \left( \frac{\min_j \nu_j}{\frac{n}{r} - 1} - 1 \right),
\end{equation}

\begin{equation}\label{cond3}
\min_j \beta_j > \frac{1}{m} \sum_{j=1}^m \beta_j - \frac{n}{rp} + \frac{m}{p} \left( 2\frac{n}{r} - 1 + \max_j \nu_j
\right),
\end{equation}
then $T_\beta$ is bounded from $\prod_{k=1}^m L^p_{m\nu_k + (m-1)\frac{n}{r}}(T_\Omega)$ to
 $L^p((T_\Omega)^m, \prod_{k=1}^m \Delta^{\nu_k - \frac{n}{r}} d V(z_k))$.
\end{thm}

Among our applications of the above result, we obtain a sufficient condition for the boundedness of the Bergman
projection in terms of the reproducing formula, which is new in this setting. More precisely, we prove the following
theorem.

\begin{thm}
Let $\nu > \frac{n}{r} - 1$ and $1 < p < \infty$. If for any $f \in L^p_\nu(T_\Omega)$ the following representation
formula holds
\begin{equation}\label{repr1}
P_\nu f(z_1) P_\nu f(z_2) = C_\beta \int_{T_\Omega}
\frac{f(z)P_\nu f(z) \Delta^{\beta - \frac{n}{r}}(\Im z)}{\Delta^{\frac{1}{2}(\frac{n}{r}+\frac{\beta}{2})}
\left(\frac{z_1 - \overline z}{i} \right)\Delta^{\frac{1}{2}(\frac{n}{r}+\frac{\beta}{2})}
\left(\frac{z_2 - \overline z}{i} \right)} d V(z)
\end{equation}
for some sufficiently large $\beta$ and all $z_1, z_2$ in $T_\Omega$, then the Bergman projection $P_\nu$ is bounded on $L^p_\nu (T_\Omega)$.
\end{thm}

In this theorem the weights $\nu$ and $\beta$ are taken real, but the result generalizes directly to the vector
weight case. The condition "$\beta$ is sufficiently large" is related to the boundedness conditions for the Bergman
kernal and determinant function. For example, a necessary condition for the boundedness of the Bergman projection
$P_\beta$ on $L^p_\nu (T_\Omega)$ is that the related Bergman kernel belongs to $L^{p^\prime, q^\prime}_\nu
(T_\Omega)$, where $1/p + 1/p^\prime = 1$, $1/q + 1/q^\prime = 1$, and this can only happen for large values of
$\beta$ for $p, q$ and $\nu$ fixed, see \cite{Seh}.

The second problem considered in this paper is the embedding relation between some generalization of the classical
Hardy spaces and the weighted mixed norm Bergman spaces in tube domains over general symmetric cones. Let
$\mathcal H^p (T_\Omega)$ denotes the holomorphic Hardy space on the tube domain, i.e. the space of holomorphic
functions $f$ on $T_\Omega$ such that
$$ \| f \|_{\mathcal H^p} = \left( \sup_{t \in \Omega} \int_{\mathbb R^n} | f(x + it)|^p d x \right)^{1/p} < \infty.$$

Following \cite{Ga}, we extend the above definition of Hardy spaces to a more general family of spaces
$\mathcal H^p_\mu (T_\Omega)$ for any locally finite and quasi-invariant measure $\mu$ supported on $\overline \Omega$.
The space $\mathcal H^p_\mu (T_\Omega)$ consists of all functions $f$ holomorphic in $T_\Omega$ satisfying
$$\| f \|_{\mathcal H^p_\mu (T_\Omega)} = \left( \sup_{t \in \Omega} \int_{T_\Omega} |f(x+i(y+t))|^p d x d \mu (y)
\right)^{1/p} < \infty.$$
In particular, if $\mu = \delta_0$, this space coincides with the classical Hardy space and if $\mu$ is the Lebesgue
measure, it coincides with the Bergman space $A^p(T_\Omega)$.

We are going to consider only those measures $\mu$ which are obtained by analytic continuation from the family $d \mu_s$ of measures
$$d \mu_s(t) = \chi_\Omega (t) \frac{\Delta_s (t)}{\Gamma_\Omega (s)} \frac{d t}{ \Delta^{n/r} (t)}, \;\; s \in
\mathbb R^n, \;\; s > g_0,$$
where $\Gamma_\Omega$ denotes the gamma function of the cone $\Omega$ defined in the next section. More precisely, in
the family $\{ \mu_s \}_{s \in \mathbb C^r}$ of tempered distributions we consider only those which are positive measures. These measures come from a characterization of Gindikin (see \cite{Ga} or Theorem VII 3.2. of \cite{FaKo})
and correspond to those $s = (s_1, \ldots, s_n) \in \mathbb C^r$ which belong to the following Wallach set
$$\Xi = \left\{ (u_1, u_1 + \frac{d}{2}{\rm sgn}\, u_1, \ldots, u_r + \frac{d}{2} ({\rm sgn}\, u_1 + \cdots +
 {\rm sgn}\, u_r)) : u_1, \ldots, u_r \geq 0 \right\}.$$

We are interested in the embedding relations between $\mathcal H^p_\mu (T_\Omega)$ and $A^{p, q}_\nu (T_\Omega)$,
$\mu = \mu_s$, $s \in \Xi$ and $\nu \in \mathbb R^r$. When $V = \mathbb R^n$ we prove the following sharp result.

\begin{thm}
Let $s \in \Xi$, $\mu = \mu_s$, $\nu \in \mathbb R^r$ and $\nu > g_0$. Then for $2 \leq p, q < \infty$ with
$\frac{q}{2} s > g_0$ we have: $\mathcal H^2_\mu (T_\Omega) \hookrightarrow A^{p,q}_\nu (T_\Omega)$ if and only if
$\frac{n}{2r} + \frac{s}{2} = \frac{\nu}{q} + \left( \frac{n}{rp}, \cdots, \frac{n}{rp} \right)$.
\end{thm}

We note that for the sufficiency part of the above theorem it suffices to prove that $\mathcal H^2_\mu (T_\Omega) \hookrightarrow A^{2, u}_{\frac{u}{2}s}(T_\Omega)$ for all $u \geq 2$ such that $\frac{u}{2}s > g_0$. This is an easy consequence of the embedding relations between Bergman spaces. The condition $\frac{u}{2}s > g_0$ shows that for $s$ fixed, $u$ should be sufficiently large and so this theorem is not applicable in all cases. Moreover, it is clear that the usual Hardy space $\mathcal H^2$ is not covered by this theorem.

\begin{thm}
Let $s \in \Xi$, $\mu = \mu_s$, $\nu \in \mathbb R^r$ and $\nu > g_0$. Then for $4 \leq p, q < \infty$ we have: $\mathcal H^2_\mu (T_\Omega) \hookrightarrow A^{p,q}_\nu (T_\Omega)$ if and only if
$\left( \frac{n}{4r}, \cdots, \frac{n}{4r} \right) + \frac{1}{4} (2s + \frac{n}{r}) = \frac{\nu}{q} + \left(
\frac{n}{rp}, \cdots, \frac{n}{rp} \right)$.
\end{thm}

Again, we note that the sufficiency part of this theorem can be reduced to the proof of the embedding
$\mathcal H^2_\mu (T_\Omega) \hookrightarrow A^{4, u}_{\frac{u}{4}(2s + \frac{n}{r})}(T_\Omega)$ for $u \geq 4$.

The necessity parts of Theorems 3 and 4 follow exactly as in Proposition 2.25 of \cite{DD} with the use of norm
identity provided in Proposition 3.1 of \cite{Ga} for $\mathcal H^2_\mu(T_\Omega)$. In order to prove sufficiency, we
heavily rely on Paley-Wiener theory in this setting. The only difference with the one-dimensional case is that one has
to deal with the beta function of the tube domain with respect to the rotated Jordan frame. That is, one needs a
version of Theorem VII 1.7. of \cite{FaKo} where the generalized determinant function is replaced by the one
corresponding to the rotated Jordan frame. This version has been obtained in the forthcoming work \cite{BGNS}.

Finally, throughout this paper $C$ or $c$ denote positive constants, not necessarily the same at different occurences;
dependence on parameters is indicated by subscripts. Given two quantities $A$ and $B$, the notation $A\lesssim B$ means
that there is an absolute constant $C$ such that $A \lesssim CB$. When both $A \lesssim B$ and $B \lesssim A$ hold we
write $A\approx B$.

\section{Preliminaries and auxiliary results}

For reader's convenience, we collect in this section some definitions and results that are used in this paper, they are
essentially contained in \cite{FaKo}.

\subsection{Symmetric cones and the generalized determinant function}

Let $\Omega$ be an irreducible open cone of rank $r$ in an $n$-dimensional vector space $V$ endowed with an inner product $(\cdot / \cdot)$ for which $\Omega$ is self-dual. Let $G(\Omega)$ be the group of transformations of $\Omega$
and $G$ its identity component. It is well known that there is a subgroup $H$ of $G$ acting simply transitively on
$\Omega$, i.e. for every $y \in \Omega$ there is a unique $g \in H$ such that $y = g \bf e$, where $\bf e$ is a fixed
element in $\Omega$.

We recall that $\Omega$ induces in $V$ a structure of Euclidean Jordan algebra with identity $\bf e$ such that
$$\overline \Omega = \{ x^2 : x \in V \}.$$
We can identify (since $\Omega$ is irreducible) the inner product $( \cdot / \cdot)$ with the one given by the trace
on $V$:
$$ (x / y) = {\rm tr}\, (xy), \;\;\; x, y \in V.$$
Let $\{c_1, \ldots, c_r\}$ be a fixed Jordan frame in $V$ and
$$V = \oplus_{1 \leq i \leq j \leq r} V_{i,j}$$
be its associated Pierce decomposition of $V$. We denote by $\Delta_1(x), \ldots, \Delta_r(x)$ the principal minors of
$x \in V$ with respect to the fixed Jordan frame $\{c_1, \ldots, c_r\}$. More precisely, $\Delta_k(x)$ is the
determinant of the projection $P_k x$ of $x$ in the Jordan subalgebra $V^{(k)} = \oplus_{1 \leq i \leq j \leq k}
V_{i,j}$. We have $\Delta = \Delta_r$ and $\Delta_k(x) > 0$, $1 \leq k \leq r$, when $x \in \Omega$. The
generalized power function on $\Omega$ is defined as
$$\Delta_s(x) = \Delta^{s_1-s_2}_1(x) \Delta^{s_2-s_3}_2(x) \cdot \Delta^{s_r}_r(x),\;\; x \in \Omega, \;
s \in \mathbb C^r.$$

Next, we recall the definition of generalized gamma function associated to $\Omega$:
$$\Gamma_\Omega(s) = \int_\Omega e^{-({\bf e}/\xi)} \Delta_s(\xi) \Delta^{-n/r}(\xi) d \xi, \;\;\;
s = (s_1, \ldots, s_r) \in \mathbb C^r.$$
This integral converges if and only if $\Re s_j > (j-1)\frac{n/r -1}{r-1} = (j-1) \frac{d}{2}$ for all $1\leq j \leq r$. In that case we have a formula:
$$\Gamma_\Omega(s) = (2\pi)^{\frac{n-r}{2}} \prod_{j=1}^r \Gamma(s_j-(j-1)\frac{d}{2}),$$
see Chapter VII of \cite{FaKo} for details. We have the following result on the Laplace transform of the generalized
power function (see Proposition VII.1.2 and Proposition VII.1.6 in \cite{FaKo}).

\begin{lem}
Let $s = (s_1, \ldots, s_r) \in \mathbb C^n$ with $\Re s_j > (j-1)\frac{d}{2}$, $j = 1, \ldots, r$. Then, for all
$y \in \Omega$ we have
$$\int_\Omega e^{-i(y/\xi)} \Delta_s (\xi) \Delta^{-n/r}(\xi) d \xi = \Gamma_\Omega (s) \Delta_s(y^{-1}) =
\Gamma_\Omega (s) [\Delta^\ast_{s^\star}(y)]^{-1}.$$
Here, $y = h \bf e$ if and only if $y^{-1} = h^{\ast -1} \bf e$ with $h \in H$ and $\Delta_j^\ast$, $j = 1, \ldots, r$
are the principal minors with respect to the rotated Jordan frame $\{c_1, \ldots, c_r \}$.
\end{lem}

The beta function of the symmetric cone $\Omega$ is defined by the following integral
$$B_\Omega(p, q) = \int_{\Omega \cap ({\bf e} - \Omega)} \Delta_{p-\frac{n}{r}}(x) \Delta_{q-\frac{n}{r}}({\bf e}- x)
d x,$$
where $p$ and $q$ are in $\mathbb C^r$. When $\Re p_j, \Re q_j > (j-1)\frac{d}{2}$ the above integral converges
absolutely and
$$B_\Omega (p, q) = \frac{\Gamma_\Omega (p) \Gamma_\Omega (q)}{\Gamma_\Omega (p+q)},$$
(see Theorem VII.1.7 in \cite{FaKo}).

Let $m$ be an element of $G_{\bf e}$, the stabilizer of $\bf e$ in $G$ such that
$$mc_j = c_{r-j+1}, \;\; j = 1, \ldots, r.$$
Then for any $y \in \Omega$ and $s \in \mathbb C^r$, $\Delta_s^\ast (y) = \Delta_s(m^{-1}y)$ (see \cite{FaKo}, page
127). Using this identity, one obtains as in the proof of Theorem VII.1.7. of \cite{FaKo} the following result (see
\cite{BGNS} for details).

\begin{lem}
Let $y \in \Omega$. The integral
$$F(y) = \int_{(y-\Omega)\cap \Omega} \Delta^\ast_{p^\star - \frac{n}{r}} (x) \Delta^\ast_{q^\star - \frac{n}{r}}(y-x) d x$$
converges if $\Re p_j, \Re q_j > (j-1) \frac{d}{2}$ for $j = 1, \ldots, r$. In this case
$$F(y) = B_\Omega(p^\star, q^\star) \Delta^\ast_{p^\star + q^\star - \frac{n}{r}} (y).$$
\end{lem}

\subsection{Bergman spaces and integrability of the Bergman kernel function}

Let us recall some estimates for the functions in the Bergman space or the projections of the functions in
$L^{p,q}_\nu (T_\Omega)$. We begin with a pointwise estimate of elements in $A^{p,q}_\nu (T_\Omega)$. The following
lemma follows from the invariance of the Bergman spaces with respect to the transformation group $G(\Omega)$ (see
\cite{DD}).

\begin{lem}
Let $1 \leq p, q < \infty$ and $\nu \in \mathbb R^r$, $\nu > g_0$. Then
\begin{equation}\label{est1}
|f(z)| \lesssim \Delta_{-\frac{\nu}{q} - \frac{n}{rp}}(\Im z) \| f \|_{A^{p,q}_\nu}, \;\;\;\; z \in T_\Omega.
\end{equation}
\end{lem}

We also need a pointwise estimate for the Bergman projection of functions in $L^{p,q}(T_\Omega)$, defined by
integral formula (\ref{bpro}), when this projection makes sense. Let us first recall the following integrability
properties for the determinant function.

\begin{lem}
Let $\alpha \in \mathbb C^r$ and $y \in \Omega$.

1) The integral
$$J_\alpha (y) =  \int_{\mathbb R^n} \left| \Delta_{-\alpha} \left( \frac{x+iy}{i} \right) \right| dx$$
converges if and only if $\Re \alpha > g_0^\ast + \frac{n}{r}$. In that case $J_\alpha (y) = C_\alpha
|\Delta_{-\alpha + n/r} (y)|$.

2) For any multi-indices $s$ and $\beta$ and $t \in \Omega$ the function $y \mapsto \Delta_\beta(y+t) \Delta_s(y)$
belongs to $L^1(\Omega, \frac{dy}{\Delta^{n/r}(y)})$ if and only if $\Re s > g_0$ and $\Re (s+\beta) < g_0^\ast$. In
that case we have
$$\int_\Omega \Delta_\beta (y) \Delta_s(y) \frac{dy}{\Delta^{n/r}(y)} = C_{\beta, s} \Delta_{s+\beta}(y).$$
\end{lem}

We refer to Corollary 2.18 and Corollary 2.19 of \cite{DD} for the proof of the above lemma. Let $\tau$ denotes the set
of all triples $(p, q, \nu)$ such that $1 \leq p, q < \infty$, $\nu > g_0$ and the function $B_\nu (\cdot, i \bf e)$
belongs to $L^{p^\prime, q^\prime}_\nu (T_\Omega)$. We have the following pointwise estimate.

\begin{lem}
Suppose $(p, q, \nu) \in \tau$. Then
\begin{equation}\label{poin}
|P_\nu f(z)| \leq \Delta_{-\frac{\nu}{q} - \frac{n}{rp}}(\Im z) \| f \|_{L^{p,q}_\nu}.
\end{equation}
\end{lem}
{\it Proof.} This is an easy consequence of the above lemma and H\"older's inequality. $\Box$

We conclude this section with a useful embedding relation between mixed norm Bergman spaces (see \cite{DD} for an
alternative proof).

\begin{lem}
Suppose $1 \leq p \leq s < \infty$, $1 \leq q \leq t < \infty$ and $\nu, \beta > g_0$. Then $A^{p,q}_\nu (T_\Omega) \hookrightarrow A^{s,t}_\beta (T_\Omega)$ if and only if $\frac{\nu}{q} + \frac{n}{rp} = \frac{\beta}{t} +
\frac{n}{rs}$.
\end{lem}

{\it Proof.} Let us suppose that $\frac{\nu}{q} + \frac{n}{rp} = \frac{\beta}{t} + \frac{n}{rs}$. We recall that there is a sequence of points $(z_{j,k} = x_{j,k} + y_k)_{j,k \in \mathbb Z}$ such that
$$\| f\|_{A^{l,m}_\mu}^m \approx \sum_k \left( \sum_j |f(z_{j,k})|^l \right)^{m/l} \Delta_{\mu + \frac{n}{r} \frac{m}{l}}(y_k),$$
see \cite{BGNS} and \cite{BBPR}. From this and embeddings between $l^p$ spaces we obtain
\begin{eqnarray*}
\| f \|^t_{A^{s,t}_\beta} & \approx &
\sum_k \left( \sum_j |f(z_{j,k})|^s \right)^{t/s} \Delta_{\beta + \frac{nt}{rs}}(y_k) \cr
& \leq & \sum_k \left( \sum_j |f(z_{j,k})|^p \right)^{t/p} \Delta_{t(\frac{\nu}{q} + \frac{n}{rp})}(y_k) \cr
& \leq & \left( \sum_k \left( \sum_j |f(z_{j,k})|^P \right)^{q/p} \Delta_{\nu + \frac{nq}{rp}}(y_k)
\right)^{t/q} \approx \| f \|^t_{A^{p,q}_\nu}. \cr
\end{eqnarray*}

For the converse, we test with functions $B_\mu (\cdot, x+iy)$ where $\mu$ is large enough and $x+iy$ is fixed in
$T_\Omega$. Now continuity of the embedding and Lemma 4 give
$$\Delta_{-t\mu + \beta + \frac{nt}{rs}}(y) \leq C \Delta_{-t\mu + \nu \frac{t}{q} + \frac{nt}{rp}}(y), \;\;\;
y \in \Omega,$$
which implies that $\frac{\nu}{q} + \frac{n}{rp} = \frac{\beta}{t} + \frac{n}{rs}$. $\Box$

As a first application of the above lemma, we see that for the proof of the sufficiency in Theorem 3 it is enough to
prove that $\mathcal H^2_\mu (T_\Omega) \hookrightarrow A^{2, u}_{\frac{u}{2}s} (T_\Omega)$ for all $u \geq 2$ such that $\frac{u}{2} > g_0$. In fact, if $p, q$ and $\nu$ satisfy the hypotheses of Theorem 3 then, by the above lemma, we have
$A_{\frac{u}{2}}^{2, u}(T_\Omega) \hookrightarrow A^{p,q}_\nu (T_\Omega)$ with $u \leq q$. Similarly, for the proof of
sufficiency in Theorem 4 it suffices to prove that $\mathcal H^2_\mu (T_\Omega) \hookrightarrow A^{4, u}_{\frac{u}{4}
(2s + \frac{n}{r})}(T_\Omega)$ for all $u \geq 4$.

\section{Bergman-type operators and multifunctional embeddings}

We denote by $\Box = \Delta(\frac{1}{i} \frac{\partial}{\partial x})$ the partial differential operator of order $r$
on $\mathbb R^n$ defined by
\begin{equation}\label{box}
\Box [e^{i(x| \xi)}] = \Delta(\xi) e^{i(x| \xi)}, \;\;\;\;\; x, \xi \in \mathbb R^n.
\end{equation}

\subsection{Multifunctional Bergman-type operators}
Now we investigate boundedness of $T_\beta$ from $\prod_{k=1}^m L^p_{m\nu_k + (m-1)\frac{n}{r}}(T_\Omega)$ to
$L^p((T_\Omega)^m, \prod_{k=1}^m \Delta^{\nu_k - \frac{n}{r}} dV(z_k))$. We apply the obtained result to
multifunctional embeddings for functions in the Bergman spaces $A^p_\nu(T_\Omega)$ where $\nu > \frac{n}{r} - 1$ and
$1 \leq p < \infty$. We begin with the following result, which is known in the case $m=1$, see \cite{BBPR}.

\begin{thm}
Let $\nu = (\nu_1, \ldots, \nu_m) \in \mathbb R^m$, $m > 1$ and $1\leq p < \infty$, $\beta = (\beta_1, \ldots, \beta_m) \in \mathbb R^m$. If the parameters satisfy the following conditions
\begin{equation}\label{usl1}
\frac{1}{m} \sum_{j=1}^m \beta_j > \frac{n}{r} - 1,
\end{equation}

\begin{equation}\label{usl2}
1 \leq p < 1 + m \left( \frac{ \min_j \nu_j}{\frac{n}{r} - 1} - 1 \right),
\end{equation}
and
\begin{equation}\label{usl3}
\min_j \beta_j > \frac{1}{m} \sum_{j=1}^m \beta_j - \frac{n}{rp} + \frac{m}{p} \left(2\frac{n}{r} - 1 +
\max_j \nu_j \right),
\end{equation}
then $T_\beta$ is bounded from $\prod_{k=1}^m L^p_{m\nu_k + (m-1)\frac{n}{r}}(T_\Omega)$ to
$L^p((T_\Omega)^m, \prod_{k=1}^m \Delta^{\nu_k - \frac{n}{r}} dV(z_k))$
\end{thm}

The idea of proof is taken from \cite{LS}.

{\it Proof.} Using H\"older inequality we obtain
\begin{eqnarray*}
|T_\beta(\overrightarrow{f}(z_1, \ldots, z_m)|^p & = &
\left| \int_{T_\Omega} \frac{ \left(\prod_{j=1}^m f_j(z) \right) \Delta^{ \frac{1}{m} \sum_{j=1}^m \beta_j}(\Im z)}
{\prod_{j=1}^m \Delta^{\frac{1}{m}(\frac{n}{r} + \beta_j)}(\frac{z_j - \overline z}{i})}
\frac{d V(z)}{\Delta^{\frac{n}{r}}(\Im z)}\right|^p \cr
& \leq & I \times J, \cr
\end{eqnarray*}
where
$$I = \int_{T_\Omega}
\frac{ \left(\prod_{j=1}^m |f_j(z)|^p \right) \Delta^{ \frac{1}{m} \sum_{j=1}^m \beta_j}(\Im z)}
{\prod_{j=1}^m |\Delta (\frac{z_j - \overline z}{i})|^{p\alpha_j}}
\frac{d V(z)}{\Delta^{\frac{n}{r}}(\Im z)},$$

$$J^{p^\prime/p} = \int_{T_\Omega}
\frac{\Delta^{ \frac{1}{m} \sum_{j=1}^m \beta_j}(\Im z)}
{\prod_{j=1}^m |\Delta(\frac{z_j - \overline z}{i})|^{p^\prime \gamma_j}}
\frac{d V(z)}{\Delta^{\frac{n}{r}}(\Im z)},$$
and
\begin{equation}\label{alga}
\alpha_j + \gamma_j = \frac{1}{m} \left(\frac{n}{r} + \beta_j \right).
\end{equation}
Let us choose $\gamma_j$ such that
$\gamma_j > \frac{1}{mp^\prime} \left( \frac{1}{m} \sum_{j=1}^m \beta_j + 2 \frac{n}{r} - 1 \right)$. Then we estimate the integral $J$ using H\"older's inequality and Lemma 4:
\begin{eqnarray*}
J^{p^\prime/p} & = & \int_{T_\Omega} \prod_{j=1}^m \left| \Delta \left( \frac{z_j-\overline z}{i}
\right)\right|^{-p^\prime\gamma_j} \Delta^{\frac{1}{m} \sum_{j=1}^m \beta_j - \frac{n}{r}} (\Im z) dV(z) \cr
& \leq & C \prod_{j=1}^m \left( \int_{T_\Omega} \left| \Delta \left( \frac{z_j - \overline z}{i} \right)
\right|^{-mp^\prime\gamma_j}\Delta^{\frac{1}{m} \sum_{j=1}^m \beta_j - \frac{n}{r}} (\Im z) dV(z) \right)^{1/m} \cr
& = & C \prod_{j=1}^m \Delta^{-p^\prime\gamma_j + \frac{1}{m^2} \sum_{j=1}^m \beta_j + \frac{n}{rm}} (\Im z_j).\cr
\end{eqnarray*}
Hence we obtained:
\begin{equation}\label{zaj}
J \leq C \prod_{j=1}^m \Delta^{-p\gamma_j + \frac{p}{m^2p^\prime} \sum_{j=1}^m \beta_j + \frac{p}{p^\prime}
\frac{n}{rm}}(\Im z_j).
\end{equation}
Using the estimate (\ref{zaj}) and Lemma 4 we finally obtain
\begin{eqnarray*}
\int_{T_\Omega} \cdots \int_{T_\Omega} \prod_{k=1}^m |T_\beta(\overrightarrow{f})(z_1, \ldots, z_m)|^p
\Delta^{\nu_k - \frac{n}{r}}(\Im z_k) dV(z_1)\cdots dV(z_m) & \leq & \cr
C\int_{T_\Omega} \left( \prod_{j=1}^m |f_j(z)|^p \right) g(z) \Delta^{\frac{1}{m} \sum_{j=1}^m \beta_j} (\Im z)
\frac{dV(z)}{\Delta^{\frac{n}{r}} (\Im z)} \cr
\end{eqnarray*}
where
\begin{eqnarray*}
g(z) & = &
\int_{T_\Omega} \cdots \int_{T_\Omega} \prod_{k=1}^m \left( \left| \Delta \left( \frac{z_k-\overline z}{i} \right) \right|^{-p\alpha_k}
\Delta^{\nu_k - \frac{n}{r} - p\gamma_j + \frac{p}{m^2p^\prime} \sum_{j=1}^k \beta_j +
\frac{np}{rmp^\prime}}(\Im z_k) \right)\cr
& &  d V(z_1) \cdots d V(z_m). \cr
\end{eqnarray*}
Note that (\ref{usl3}) implies $p\alpha_k > \nu_k - p\gamma_k + \frac{p}{m^2p^\prime} \sum_{k=1}^m \beta_k +
\frac{pn}{rmp^\prime} + 2\frac{n}{r} - 1$. Thus, if we finally choose $\alpha_j$ and $\gamma_j$ such that
(\ref{alga}) holds and, for every $j = 1, \ldots, m$, we have
\begin{eqnarray*}\label{final}
\frac{1}{mp^\prime} \left( \frac{1}{m} \sum_{j=1}^m \beta_j + 2\frac{n}{r} - 1 \right) < \gamma_j & & \cr
< \min \left\{ \frac{1}{m} (\frac{n}{r} + \beta_j), \frac{\min_j \nu_j - \frac{n}{r} + 1}{p} +
\frac{\frac{1}{m} \sum_{j=1}^m \beta_j + \frac{n}{r}}{mp^\prime} \right\}, & & \cr
\end{eqnarray*}
then an application of Lemma 4 gives estimate
$$g(z) \leq C \Delta^{\sum_{k=1}^m \nu_k + m \frac{n}{r} - p \sum_{k=1}^m (\alpha_k + \beta_k) +
\frac{p}{mp^\prime} \sum_{k=1}^m \beta_k + \frac{pn}{rp^\prime}} (\Im z).$$
Finally, using H\"older's inequality we obtain
\begin{eqnarray*}
\int_{T_\Omega} \cdots \int_{T_\Omega} \prod_{k=1}^m |T_\beta (\overrightarrow{f})(z_1, \ldots, z_m)|^p
\Delta^{\nu_k - \frac{n}{r}}(\Im z_k) d V(z_1) \cdots d V(z_m) & \leq & \cr
C\int_{T_\Omega} \left( \prod_{j=1}^m |f_j(z)|^p \right) \Delta^{\sum_{k=1}^m \nu_k + (m-1) \frac{n}{r}}
(\Im z) \frac{d V(z)}{\Delta^{\frac{n}{r}} (\Im z)} & \leq & \cr
C \left( \prod_{j=1}^m \int_{T_\Omega} |f_j(z)|^{mp} \Delta^{m\nu_j + (m-1)\frac{n}{r}}(\Im z) \frac{d V(z)}{
\Delta^{\frac{n}{r}}(\Im z)} \right)^{1/m} & < & \infty. \Box \cr
\end{eqnarray*}

An analogue of the following lemma in the setting of the unit ball in $\mathbb C^n$ is contained in \cite{LS}.
Note that the case $m=1$ is obvious.

\begin{lem}
Let $\nu_k > \frac{n}{r}-1$, $k=1, \ldots, m$ and $1 \leq p < \infty$. Then there is a constant $C > 0$ such that
\begin{equation}\label{lem7}
\int_{T_\Omega} \prod_{k=1}^m |f_k(z)|^p \Delta^{(m-1)\frac{n}{r} + \sum_{k=1}^m \nu_k - \frac{n}{r}} (\Im z)
d V(z) \leq C \prod_{k=1}^m \| f_k \|^p_{A^p_{\nu_k}}.
\end{equation}
\end{lem}

{\it Proof.} By Lemma 6 we have $A^{p/m}_{\frac{1}{m}\sum_{k=1}^m \nu_k} (T_\Omega) \hookrightarrow
A^p_{(m-1)\frac{n}{r} + \sum_{k=1}^m \nu_k} (T_\Omega)$. Thus, to prove the lemma, we only need to check that for
$f_j \in A^p_{\nu_j}(T_\Omega)$, $j=1, \ldots, m$, the product $f_1 \cdots f_m$ is in
$A^{p/m}_{\frac{1}{m}\sum_{k=1}^m \nu_k} (T_\Omega)$ with the appropriate norm estimate. An application of H\"older's
inequality
$$\int_{T_\Omega} \prod_{k=1}^m |f_k(z)|^p \Delta^{\frac{1}{m}\sum_{k=1}^m \nu_k - \frac{n}{r}} (\Im z) d V(z)  \leq
\prod_{k=1}^m \left( \int_{T_\Omega} |f_k(z)|^p \Delta^{\nu_k - \frac{n}{r}} (\Im z) d V(z) \right)^{1/m}$$
finishes the proof since the last expression is equal to $\prod_{k=1}^m \| f_k \|^{p/m}_{A^p_{\nu_k}}$. $\Box$

A complete analogue of the following multifunctional result in the setting of the unit ball in $\mathbb C^n$ can be
found in \cite{LS}.

\begin{thm}
Let $\nu_k > \frac{n}{r}$ for $1 \leq k \leq m$, $m>1$. Let $1 \leq p < \infty$ and suppose that $\beta_j$ are
sufficiently large so that for any sequence $(z_j)_{j=1}^m$ in $T_\Omega$ the following representation holds for $f_1, \ldots, f_m \in \mathcal H(T_\Omega)$
\begin{equation}\label{repr2}
f_1(z_1) \cdots f_m(z_m) = C_{m, \beta} \int_{T_\Omega}
\frac{\prod_{j=1}^m f_j(z)\Delta^{\frac{1}{m} \sum_{j=1}^m} (\Im z)}{\prod_{j=1}^m \Delta^{\frac{1}{m}
(\frac{n}{r} + \beta_j)} (\frac{z_j-\overline z}{i})} \frac{ d V(z)}{\Delta^{n}{r} (\Im z)}.
\end{equation}
Assuming none of the functions $f_k$ is identically zero, the following statements are equivalent.

1) There is a constant $C > 0$ such that
\begin{equation}\label{bdd}
\int_{T_\Omega} \prod_{k=1}^m |f_k(z)|^p \Delta^{(m-1)\frac{n}{r} + \sum_{k=1}^m \nu_k} (\Im z) \frac{ d V(z)}{
\Delta^{\frac{n}{r}}(\Im z)} \leq C < \infty.
\end{equation}

2) $f_k \in A^p_{\nu_k}(T_\Omega)$ for all $k=1, \ldots, m$.
\end{thm}

{\it Proof.} We have already seen that $2) \Rightarrow 1)$ independently of the representation formula (\ref{repr2}). Let us prove implication $1) \Rightarrow 2)$ assuming (\ref{repr2}). Since the functions $f_j$ are not identically
zero, condition
$$\int_{T_\Omega} \cdots \int_{T_\Omega} \prod_{k=1}^m (|f_k(z_k)||^p \Delta^{\nu_k - \frac{n}{r}} (\Im z_k)
d V(z_1) \cdots d V(z_m) < \infty$$
implies $f_k \in A^p_{\nu_k}(T_\Omega)$ for all $k = 1, \ldots, m$. Now, using the representation (\ref{repr2}) we
obtain
\begin{eqnarray*}
K & = & \int_{T_\Omega} \cdots \int_{T_\Omega} \prod_{k=1}^m (|f_k(z_k)|^p \Delta^{\nu_k - \frac{n}{r}} (\Im z_k))
d V(z_1) \cdots d V(z_m) \cr
& = & \int_{T_\Omega} \cdots \int_{T_\Omega} |T_\beta (\overrightarrow{f}(z_1, \ldots, z_m)|^p
\left( \prod_{k=1}^m \Delta^{\nu_k - \frac{n}{r}} (\Im z_k) \right) d V(z_1) \cdots d V(z_m), \cr
\end{eqnarray*}
where $\overrightarrow{f} = (f_1, \ldots, f_m)$. The proof of Theorem 5 gives
$$ K \leq C \int_{T_\Omega} \cdots \int_{T_\Omega} \prod_{k=1}^m (|f_k(z_k)||^p \Delta^{\nu_k - \frac{n}{r}} (\Im z_k)
d V(z_1) \cdots d V(z_m) < \infty. \;\;\; \Box $$

We write $(\nu, p) \in \sigma$ if $\nu \in \mathbb R$, $1 \leq p < \infty$, $\nu > \frac{n}{r} - 1$ and
$\Delta^{-(\nu + \frac{n}{r})}(\frac{z-i{\bf e}}{i}) \in L^{p^\prime}_\nu (T_\Omega)$. Let us define, for
$f_k \in L^p_{\nu_k}$, the following operators:
\begin{equation}\label{sbk}
S_{\beta, k}(\overrightarrow{f})(\overrightarrow{z}) = \int_{T_\Omega}
\frac{f_k(z) \prod_{j \not= k} P_{\nu_j}f_j(z) \Delta^{\frac{1}{m} \sum_{j=1}^m \beta_j} (\Im z)}{\prod_{j=1}^m
\Delta^{\frac{1}{m}(\frac{n}{r} + \beta_j)}(\frac{z_j - \overline z}{i})} \frac{ d V(z)}{\Delta^{\frac{n}{r}} (\Im z)},
\end{equation}
and
\begin{equation}\label{sb}
S_\beta = \sum_{k=1}^m S_{\beta, k}.
\end{equation}

\begin{thm}
Suppose $(\nu_k, p) \in \sigma$ for $k=1, \ldots, m$. If the parameters satisfy conditions (\ref{cond1}), (\ref{cond2}), and (\ref{cond3}), then the operators $S_{\beta_k}$ and $S_\beta$ are bounded from $\prod_{j=1}^m L^p_{\nu_j}
(T_\Omega)$ to $L^p((T_\Omega)^m, \prod_{j=1}^m  \Delta^{\nu_j - \frac{n}{r}} (\Im z_j) d V(z_j))$.
\end{thm}

{\it Proof.} Clearly we only need to prove the result for $S_{\beta_k}$ for fixed $k$. An inspection of the proof of
Theorem 5 and Lemma 5 give
\begin{eqnarray*}
\int_{T_\Omega} \cdot \int_{T_\Omega} |S_{\beta, k}(\overrightarrow{f}(\overrightarrow{z}|^p
\Delta^{\nu_j - \frac{n}{r}}(\Im z_j) d V(z_1) \cdots d V(z_m) & \leq & \cr
C \int_{T_\Omega} |f_k(z)|^p \left( \prod_{j \not= k}^m |P_{\nu_j} f_j(z)|^p \right)
\Delta^{\sum_{k=1}^m \nu_k + (m-1)\frac{n}{r}}(\Im z) \frac{ d V(z)}{\Delta^{\frac{n}{r}} (\Im z)} & \leq & \cr
C \prod_{j\not= k} \| f_j \|_{L^p_{\nu_j}} \int_{T_\Omega} |f_k(z)|^p \Delta^{\nu_k - \frac{n}{r}} (\Im z) d V(z)
& \leq & C \prod_{j\not= k} \| f_j \|_{L^p_{\nu_j}}, \cr
\end{eqnarray*}
and the proof is complete. $\Box$

As a consequence we have the following result.

\begin{thm}
Suppose $(\nu_k, p) \in \sigma$ for $k=1, \ldots, m$. Suppose also that, for $\beta_j$ large enough, the
following representation
\begin{equation}\label{prod}
\prod_{k=1}^m P_{\nu_k} f_k(z_k) = C_{m, \beta} \int_{T_\Omega}
\frac{f_k(z) \prod_{j\not= k}^m P_{\nu_j} f_j(z) \Delta^{\frac{1}{m} \sum_{j=1}^m \beta_j}(\Im z)}{\prod_{j=1}^m
\Delta^{\frac{1}{m} ( \frac{n}{r} + \beta_j)} (\frac{z_j- \overline z}{i})} \frac{d V(z)}{\Delta^{\frac{n}{r}}
(\Im z)}
\end{equation}
holds for any sequence $(z_j)_{j=1}^m$ in $T_\Omega$ and any $f_k \in L^p_{\nu_k}(T_\Omega)$, $1 \leq k \leq m$.
Then $P_{\nu_k} f_k \in L^p_{\nu_k}(T_\Omega)$, $1 \leq k \leq m$.
\end{thm}

We also have the following corollary which gives a sufficient condition for boundedness of the Bergman projection.

\begin{cor}
Let $(\nu, p) \in \sigma$. If the following representation
\begin{equation}\label{rep}
P_\nu f(z_1) P_\nu f(z_2) = C_\beta \int_{T_\Omega}
\frac{f(z)P_\nu f(z)}{\Delta^{\frac{1}{2}(\frac{n}{r} + \frac{\beta}{2})} (\frac{z_1 - \overline z}{i})
\Delta^{\frac{1}{2}(\frac{n}{r} + \frac{\beta}{2})} (\frac{z_2 - \overline z}{i})} \Delta^{\beta - \frac{n}{r}}
(\Im z) d V(z)
\end{equation}
holds for all $z_1, z_2 \in T_\Omega$ and $f \in L^p_\nu(T_\Omega)$, where $\beta$ is large enough, then $P_\nu$ is
bounded on $L^p_\nu (T_\Omega)$.
\end{cor}

{\it Proof.} Using Lemma 5 we clearly have
\begin{eqnarray*}
\int_{T_\Omega} |f(z)|^p |P_\nu f(z)|^p \Delta^{2\nu} (\Im z) d V(z) & \leq & C \| f ||^p_{L^p_\nu}
\int_{T_\Omega} |f(z)|^p \Delta^{\nu - \frac{n}{r}} (\Im z) d V(z) \cr
& = & C \| f \|^{2p}_{L^p_\nu}.
\end{eqnarray*}
Now, following the proof of Theorem 7 we obtain
\begin{eqnarray*}
\| P_\nu f \|^{2p}_{L^p_\nu} & = & \int_{T_\Omega} \int_{T_\Omega} |P_\nu f(z_1)|^p |P_\nu f(z_2)|^p
\Delta^{\nu - \frac{n}{r}} (\Im z_1)\Delta^{\nu - \frac{n}{r}} (\Im z_2) d V(z_1) d V(z_2) \cr
& \leq & C \int_{T_\Omega} |f(z)|^p |P_\nu f(z)|^p \Delta^{2\nu} (\Im z) d V(z) \cr
& \leq & C \| f \|^{2p}_{L^p_\nu}. \;\;\;\; \Box \cr
\end{eqnarray*}

\subsection{Multifunctional inequalities involving Bergman projection or the box operator} Next we derive
multifunctional inequalities involving the Bergman projection or the box operator. As a preparation, we first prove
the following elementary proposition.

\begin{prop}
Let $(\nu, p) \in \sigma$. If $P_\nu$ is bounded on $L^p_\nu(T_\Omega)$, then $P_\nu$ is bounded from
$L^p_\nu (T_\Omega)$ to $L^{kp}_{k\nu + (k-1)\frac{n}{r}}(T_\Omega)$ for any $k \in \mathbb N$.
\end{prop}

{\it Proof.} Suppose $P_\nu$ is bounded on $L^p_\nu (T_\Omega)$. Then using Lemma 5 we obtain, for any $f \in
L^p_\nu (T_\Omega)$:
\begin{eqnarray*}
\int_{T_\Omega} |P_\nu f(z)|^{kp} \Delta^{k\nu + (k-2)\frac{n}{r}} (\Im z) d V(z) & = & \int_{T_\Omega}
\left( |P_\nu f(z)|^p \Delta^{\nu + \frac{n}{r}}(\Im z) \right)^{k-1} |P_\nu f(z)|^p \cr
& & \Delta^{\nu - \frac{n}{r}} (\Im z) d V(z) \cr
& \leq & C \| f \|^{(k-1)p}_{L^p_\nu} \int_{T_\Omega} |P_\nu f(z)|^p \Delta^{\nu - \frac{n}{r}} (\Im z) d V(z) \cr
& \leq & C \| f \|^{kp}_{L^p_\nu}. \;\;\;\;\; \Box
\end{eqnarray*}

\begin{prop}
Let $(\nu_k, p) \in \sigma$ for $1\leq k \leq m$. Suppose $P_{\nu_k}$ is bounded on $L^p_{\nu_k}(T_\Omega)$ for all
$k = 1, \ldots, m$. Then for any $l \in \mathbb N$ we have
$$\int_{T_\Omega} \prod_{k=1}^m [|P_{\nu_k}| f_k (z)|^{lp} \Delta^{l\nu_k + l \frac{n}{r}}(\Im z)]
\frac{d V(z)}{\Delta^{2\frac{n}{r}}(\Im z)} \leq C \prod_{k=1}^m \| f_k \|^{lp}_{L^p_{\nu_k}}.$$
\end{prop}

{\it Proof.} Using the above proposition, H\"older's inequality and Lemma 5 we obtain
\begin{eqnarray*}
\int_{T_\Omega} \prod_{k=1}^m [|P_{\nu_k}| f_k (z)|^{kp} \Delta^{l\nu_k + l \frac{n}{r}}(\Im z)]
\frac{d V(z)}{\Delta^{2\frac{n}{r}}(\Im z)} & \leq & \cr
C \prod_{k=1}^m \| f_k \|^{(l-1)p}{L^p_{\nu_k}} \int_{T_\Omega}\prod_{k=1}^m [|P_{\nu_k}| f_k (z)|^p \Delta^{\nu_k +  \frac{n}{r}}(\Im z)] \frac{d V(z)}{\Delta^{2\frac{n}{r}}(\Im z)} & \leq & \cr
C \prod_{k=1}^m \| f_k \|^{(l-1)p}_{L^p_{\nu_k}} \prod_{k=1}^m \left( \int_{T_\Omega} |P_{\nu_k} f_k(z)|^{mp}
\Delta^{m\nu_k + (m-2)\frac{n}{r}} (\Im z) d V(z) \right)^{1/m} & \leq &  \cr
C \prod_{k=1}^m \| f_k \|^{lp}_{L^p_{\nu_k}} & & \Box
\end{eqnarray*}

It is well-known that the operator $\Box$ satisfies the following boundedness estimate
\begin{equation}\label{box}
\| \Box f \|_{A^p_{\nu + p}} \leq C \| f \|_{A^p_\nu},
\end{equation}
see \cite{BBPR}. It follows, using H\"older's inequality, that for $1 \leq p < \infty$ and $q<p$
\begin{equation}\label{boxpq}
\int_{T_\Omega} | \Box f(z) |^q |f(z)|^{p-q} \Delta^{\nu + q - \frac{n}{r}} (\Im z) d V(z) \leq C \| f \|^p_{A^p_\nu}.
\end{equation}
Our goal is to obtain a multifunctional version of the above estimate. To this end, we introduce the following operator, which we still denote by $\Box$, defined for pointwise products of holomorphic functions:
$$\Box (f_1 \cdots f_m) = \sum_{j=1}^m f_1 \cdots f_{j-1} (\Box f_j) f_{j+1} \cdots f_m.$$
We note that the $\Box$ inside the sum is the usual $\Box$ as defined at the beginning of this section. The next theorem generalizes (\ref{boxpq}), this idea appeared in \cite{LS}.

\begin{thm}
Let $\nu > \frac{n}{r} -1$, $1 \leq q \leq p < \infty$. Then there exists $C > 0$ such that
\begin{equation}\label{boxes}
\int_{T_\Omega} |\Box (f_1 \cdots f_m)|^q \prod_{j=1}^m |f_j(z)|^{p-q} \Delta^{m(\nu + \frac{n}{r})+q} (\Im z)
\frac{d V(z)}{\Delta^{2\frac{n}{r}} (\Im z)} \leq C m^q \prod_{j=1}^m \| f_j \|^p_{A_\nu^p}.
\end{equation}
\end{thm}

{\it Proof.} Using Minkowski's inequality, the pointwise estimate for functions in $A^p_\nu (T_\Omega)$ and the
estimate (\ref{boxpq}) we obtain
\begin{eqnarray*}
\int_{T_\Omega}|\Box (f_1 \cdots f_m)|^q \prod_{j=1}^m |f_j(z)|^{p-q} \Delta^{m(\nu + \frac{n}{r})+q} (\Im z)
\frac{d V(z)}{\Delta^{2\frac{n}{r}} (\Im z)} & \leq & \cr
C \left( \sum_{j=1}^m \left( \int_{T_\Omega} \prod_{k\not= j}^m |f_k(z)|^q |\Box f_j(z)|^q \prod_{k=1}^m
|f_k(z)|^{p-q} \Delta^{m(\nu + \frac{n}{r})+q} (\Im z) \frac{d V(z)}{\Delta^{2\frac{n}{r}}(\Im z)} \right)^{1/q}
\right)^q & \leq & \cr
C \left( \sum_{j=1}^m \left( \int_{T_\Omega} \left( \prod_{k\not= j}^m |f_k(z)|^p \Delta^{\nu + \frac{n}{r}}
(\Im z) \right) |\Box f_j(z)|^q |f_j(z)|^{p-q} \Delta^{\nu - \frac{n}{r} + q} (\Im z) d V(z) \right)^{1/q} \right)^q
& \leq & \cr
C \left( \sum_{j=1}^m \left( \prod_{k\not= j}^m \| f_k \|_{A^p_\nu}^{p/q} \right) \left( \int_{T_\Omega}
|\Box f_j(z)|^q |f_j(z)|^{p-q} \Delta^{\nu - \frac{n}{r} + q} (\Im z) d V(z) \right)^{1/q} \right)^q & \leq & \cr
Cm^q\prod_{k=1}^m \| f_k \|^p_{A^p_\mu}. & & \cr
\end{eqnarray*}
$\Box$

\section{Paley-Wiener representation and embeddings}

We make use of Paley-Wiener theory in this section to prove Theorem 3 and Theorem 4. From now on we fix a measure
$\mu = \mu_s$, where $s \in \Xi$. We recall that $\mathcal H^2_\mu(T_\Omega)$ is a Hilbert space and we use notation
from \cite{Ga}:
$$L^2_{s^\star} (\Omega) = L^2(\Omega, \Delta^\ast_{s^\star}(2\xi)d \xi) = L^2(\Omega, \Delta_s((2\xi)^{-1})^{-1}
d \xi).$$
The following Paley-Wiener characterization of functions in $\mathcal H^2_\mu (T_\Omega)$ has been obtained in
\cite{Ga}.

\begin{thm}
For every $F \in \mathcal H^2_\mu (T_\Omega)$ there is an $f \in L^2_{s^\star}(\Omega)$ such that
$$F(z) = \frac{1}{(2\pi)^{\frac{n}{2}}} \int_\Omega e^{i(x/\xi)} f(\xi) \Delta^\ast_{s^\star} (2\xi) d \xi, \;\;\;
z \in T_\Omega.$$
Conversely, if $f \in L^2_{s^\star}(\Omega)$ then the above integral converges absolutely to a function $F \in
\mathcal H^2_\mu(\Omega)$. Moreover, $\| F \|_{\mathcal H^2_\mu} = \| f \|_{L^2_{s^\star}}$.
\end{thm}

As remarked in the introductory section, we only need to show the following result in proving Theorem 3.

\begin{thm}
Let $s \in \Xi$, $\mu = \mu_s$. For all $2 \leq q < \infty$ such that $\frac{q}{2} s > g_0$ we have
$$\mathcal H^2_\mu (T_\Omega)\hookrightarrow A^{2,q}_{\frac{q}{2} s} (T_\Omega).$$
\end{thm}

{\it Proof.} Let $F \in \mathcal H^2_\mu (T_\Omega)$, by Theorem 10 there is an $f$ in $L^2_{s^\star}(\Omega)$ such
that
$$F(z) = C_n \int_\Omega e^{i(x/\xi)} f(\xi) \Delta^\ast_{s^\star} (2\xi) d \xi, \;\;\;\;\; z \in T_\Omega.$$
It follows from Plancherel's theorem that
$$\int_{\mathbb R^n} |F(x+iy)|^2 dx = C \int_\Omega e^{-2(y/\xi)} |f(\xi)|^2 \Delta^\ast_{2s^\star}(\xi) d \xi.$$
Integrating the $q/2$-power of the left hand side of the above equality with respect to the measure
$\Delta_{\frac{q}{2}s}(y)\Delta^{-\frac{n}{r}}(y) d y$ and using Minkowski's inequality for integrals and Lemma 1
we obtain
\begin{eqnarray*}
I & = & C \int_\Omega \left( \int_\Omega e^{-2(y/\xi)} |f(\xi)|^2 \Delta^\ast_{2s^\star}(\xi) d \xi \right)^{q/2}
\Delta_{\frac{q}{2}s}(y) \frac{dy}{\Delta^{n/r}(y)} \cr
& \leq & C \left( \int_\Omega \left( \int_\Omega e^{-q(y/\xi)} \Delta_{\frac{q}{2}s}(y) \frac{dy}{\Delta^{n/r}(y)}
\right)^{2/q} |f(\xi)|^2 \Delta^\ast_{2s^\star}(\xi) d \xi \right)^{q/2} \cr
& = & C \left( \int_\Omega \Delta^\ast_{-s^\star}(\xi) |f(\xi)|^2 \Delta^\ast_{2s^\star}(\xi) d \xi \right)^{q/2} =
C \| f \|^q_{L^q_{s^\star}}, \cr
\end{eqnarray*}
where
$$I = \int_\Omega \left( \int_{\mathbb R^n} |F(x+iy)|^2 dx \right)^{\frac{q}{2}} \Delta_{\frac{q}{2}s}(y)
\frac{dy}{\Delta^{n/r}(y)} = \| F \|_{A^{2,q}_{\frac{q}{2}s}}^q. \;\;\;\;\;\; \Box$$

For our last result we need the following Paley-Wiener construction of functions in the Bergman space $A^{2,q}_\nu$.

\begin{lem}
Let $2 \leq q < \infty$ and $\nu \in \mathbb R^r$, $\nu > g_0$. If $f$ is in the space
$L^2_{2(1-\frac{1}{q})\nu^\star}(\Omega) = L^2(\Omega, \Delta^\ast_{2(1-\frac{1}{q})\nu^\star}(2\xi) d \xi)$, then
the function $F$ defined by
\begin{equation}\label{pw}
F(z) = \frac{1}{(2\pi)^{\frac{n}{2}}} \int_\Omega e^{i(x/\xi)} f(\xi) \Delta^\ast_{s^\star}(2\xi) d \xi, \;\;\;\;\;
z \in T_\Omega
\end{equation}
belongs to $A^{2,q}_\nu (T_\Omega)$.
\end{lem}

{\it Proof.} The estimation of the $L^{2,q}_\nu$-norm of the integral in (\ref{pw}) proceeds exactly as in the
previous theorem and one obtains
$$ \| F \|_{A^{2,q}_\nu} \leq C \| f \|_{L^2_{2(1-\frac{1}{q})\nu^\star}}.$$
Thus, we only have to prove that for any $f \in L^2_{2(1-\frac{1}{q})\nu^\star}(\Omega)$ the integral in
(\ref{pw}) converges absolutely to a holomorphic function $F(z)$ on $T_\Omega$. It suffices to prove this at the
point $z = i\bf e$. Using H\"older's inequality and Lemma 1 we obtain
\begin{eqnarray*}
\int_\Omega e^{-({\bf e}/\xi)} |f(\xi)| \Delta^\ast_{s^\star}(2\xi) d \xi & \leq &
\| f \|_{L^2_{2(1-\frac{1}{q})\nu^\star}} \left( \int_\Omega e^{-2({\bf e}/\xi)} \Delta^\ast_{\frac{2}{q}\nu^\star}
(2\xi) d \xi \right)^{1/2} \cr
& = & \| f \|_{L^2_{2(1-\frac{1}{q})\nu^\star}} 2^{-\frac{n}{2}} \Gamma_\Omega \left( \frac{2}{q} \nu^\star +
\frac{n}{r} \right)^{1/2}, \cr
\end{eqnarray*}
and this is clearly finite. $\Box$

We now give a proof of the following result, which, as remarked in the first section, implies Theorem 4.

\begin{thm}
Let $s \in \Xi$, $\mu = \mu_s$. For all $4 \leq q < \infty$ we have
$$\mathcal H^2_\mu (T_\Omega) \hookrightarrow A^{4,q}_{\frac{q}{4}(2s + \frac{n}{r})}(T_\Omega).$$
\end{thm}

{\it Proof.} Given $F$ in $\mathcal H^2_\mu (T_\Omega)$ we need to show that $F^2$ belongs to $A^{2,q/2}_{\frac{q}{4}
(2s + \frac{n}{r})} (T_\Omega)$. By Theorem 10 there exists $f \in L^2_{s^\star}(\Omega)$ such that
$$F(z) = C_n \int_\Omega e^{i(x/\xi)} f(\xi) \Delta^\ast_{s^\star} (2\xi) d \xi, \;\;\;\;\; z \in T_\Omega.$$
Using this Paley-Wiener representation we get
\begin{eqnarray*}
F^2(z) & = & C^2_n \int_{\Omega \times \Omega} e^{i(x/\xi + t)} f(\xi)f(t) \Delta^\ast_{s^\star}(2\xi)
\Delta^\ast_{s^\star}(2t) d\xi d t \cr
& = & C_n^2 \int_\Omega \int_{\Omega \cap (u-\Omega)} e^{i(x/u)} f(u-\xi)f(\xi) \Delta^\ast_{s^\star}
(2(u-\xi))\Delta^\ast_{s^\star}(2\xi) d\xi d u \cr
& =& C_n^2 \int_\Omega e^{i(x/u)} g(u) du, \cr
\end{eqnarray*}
where
$$g(u) = \int_{\Omega \cap (u-\Omega)} f(u-\xi)f(\xi) \Delta^\ast_{s^\star}(2(u-\xi)) \Delta^\ast_{s^\star}
(2\xi) d\xi.$$
By Lemma 8 it suffices to prove that $g(u)\Delta^\ast_{-\frac{q}{4}(2s^\star + \frac{n}{r})}(u)$ is in
$L^2_{(\frac{q}{2} - 1)(2s^\star + \frac{n}{r})}(\Omega)$, or, equivalently, that $g$ is in
$L^2_{-(2s^\star + \frac{n}{r})}(\Omega)$. We start with a pointwise estimate of $g(u)$. Using H\"older's
inequality and Lemma 2 we obtain
\begin{eqnarray*}
|g(u)|^2 & \leq & \left( \int_{\Omega \cap (u-\Omega)} |f(u-\xi)| |f(\xi)| \Delta^\ast_{s^\star}(2(u-\xi))
\Delta^\ast_{s^\star}(2\xi) d \xi \right)^2 \cr
& \leq & \left( \int_{\Omega \cap (u-\Omega)} |f(u-\xi)|^2 |f(\xi)|^2\Delta^\ast_{s^\star}(2(u-\xi))
\Delta^\ast_{s^\star}(2\xi)\right) \cr
& & \times \left( \int_{\Omega\cap (u-\Omega)} \Delta^\ast_{s^\star}(2(u-\xi)) \Delta^\ast_{s^\star}(2\xi) d\xi \right) \cr
& = & C \Delta^\ast_{2s^\star + \frac{n}{r}}(u) \left( \int_{\Omega \cap (u-\Omega)} |f(u-\xi)|^2 |f(\xi)|^2
\Delta^\ast_{s^\star}(2(u-\xi)) \Delta^\ast_{s^\star}(2\xi) d\xi \right). \cr
\end{eqnarray*}
It easily follows that
\begin{eqnarray*}
\int_\Omega \frac{|g(u)|^2 du}{\Delta^\ast_{2s^\star + \frac{n}{r}}(u)} & \leq & C \int_\Omega
\int_{\Omega \cap (u-\Omega)} |f(u-\xi)|^2 |f(\xi)|^2 \Delta^\ast_{s^\star}(2(u-\xi)) \Delta^\ast_{s^\star}(2\xi)
d\xi du \cr
& = & C \| f \|^4_{L^2_{s^\star}} = C \| F \|^4_{\mathcal H^2_\mu} \cr
\end{eqnarray*}
and the proof is complete. $\Box$

\noindent Department of Mathematics, University of Belgrade, Studentski trg 16, 11000 Belgrade, Serbia

\noindent {\it E-mail}: {\tt arsenovic@matf.bg.ac.rs}
\smallskip

Romi F. Shamoyan

\noindent Bryansk University, Bryansk, Russia

\noindent {\it E-mail}: {\tt rshamoyan@yahoo.com}


\begin{thebibliography}{99}

\bibitem{BeBo} D. B\'ekoll\'e, A. Bonami, {\it Analysis on tube domains over light cones: some extensions of recent
results}, Actes des rencontres d'analyse complexe. Poitiers (1999).

\bibitem{BBGR} D. B\'ekoll\'e, A. Bonami, G. Garrig\'{o}s, F. Ricci, {\it Littlewood-Paley decompositions related to
symmetric cones and Bergman projections in tube domains}, Proc. London Math. Soc. 89 (2004), 317-360.

\bibitem{BBGRS} D. B\'ekoll\'e, A. Bonami, G. Garrig\'{o}s, F. Ricci, B. F. Sehba, {\it Analytic Besov spaces and
Hardy-type inequalities in tube domanis over symmetric cones}, arXiv:0902.2928.

\bibitem{BBPR} D. B\'ekoll\'e, A. Bonami, M. Peloso, F. Ricci, {\it Besov spaces in tube domains over symmetric cones}

\bibitem{BGNS} A. Bonami, S. Grellier, C. Nana, B. F. Sehba, {\it Schatten classes of Hankel operators on tube
domains}, preprint.

\bibitem{DD} D. Debertol, {\it Besov spaces and boundedness of weighted Bergman projections over symmetric tube
domains}, Dottorato di Ricerca in Matematica, Universita di Genova, Politecnico di Torino (April 2003).

\bibitem{FaKo} Faraut, J. and Koranyi, A., {\it Analysis on symmetric cones}, Oxford Mathematical Monographs, Oxford SciencePublications. The Clarendon Press, Oxford University Press, New York, 1994. xii + 382 pp. ISBN:0-19-853477-9.

\bibitem{Ga} G. Garrig\'os, {\it M\"obius invariance of analytic Besov spaces on tubes over cones}, preprint.

\bibitem{GS} G. Garrig\'{o}s, A. Seeger, {\it Plate decompositions for cone multipliers}, Proc. Harmonic Analysis and
its Applications, Sapporo 2005, 13-28, Report 103.

\bibitem{LS} S. Li, R. Shamoyan, {\it On some extensions of theorems on atomic decomposition of Bergman and Bloch
spaces in the unit ball and related problems}, Comp. Var. and Ell. Equ. Vol. 52 (2009) no. 12, 1151-1162.

\bibitem{Seh} B. F. Sehba, {\it Bergman type operators in tubular domains over symmetric cones}, Proc. Edin. Math. Soc.
(Series 2), 52, (2009), pp 529-544.

\end{thebibliography}
\end{document}